\def\zbar{\overline{z}}

\def\IH{{\mathbb H}^3}

\def\IR{{\mathbb R}}

\def\ID{{\mathbb D}}

\def\IS{{\mathbb S}}
\def\IZ{{\mathbb Z}}
\def\IH{{\mathbb H}}
\def\IC{\mathbb C}

\def\zbar{{\overline{z}}}
\def\hbar{{\overline{h}}}

\documentclass{article}

\usepackage{amsthm, amssymb}
\usepackage{amsmath}
\usepackage{amsfonts}
\usepackage{epsfig,multicol}
\usepackage{pstricks}
\usepackage{graphicx}

\title{Random Lattices,  Punctured Tori and the Teichm\"uller distribution.}

\author{Gaven Martin  \thanks{Research supported in
part by grants from the N.Z.
Marsden Fund.   \newline \newline AMS
(1991) Classification.
Primary 30C60, 30F40, 30D50, 20H10, 22E40, 53A35, 57N13, 57M60}  }

\date{}

\begin{document}

\maketitle
\newtheorem{theorem}{Theorem}[section]    
                                           
\newtheorem{lemma}[theorem]{Lemma}         
\newtheorem{corollary}[theorem]{Corollary} 
\newtheorem{remark}[theorem]{Remark}       
\newtheorem{definition}[theorem]{Definition}
\newtheorem{conjecture}[theorem]{Conjecture}
\newtheorem{proposition}[theorem]{Proposition}
\newtheorem{example}[theorem]{Example}

\newcommand{\param}{(\gamma,\beta,\beta')}
\newcommand{\parfour}{(\gamma,\beta,-4)}

\numberwithin{equation}{section}

%    Absolute value notation
\newcommand{\abs}[1]{\lvert#1\rvert}

\renewcommand{\theequation}{\thetheorem} 
                    
\makeatletter
\let \c@equation=\c@theorem

\begin{abstract}  The moduli space of lattices of $\IC$ is a Riemann surface of finite hyperbolic area with the square lattice as an origin.   We select a lattice from the induced uniform distribution and calculate the statistics of the Teichm\"uller distance to the origin.  This in turn identifies distribution of the distance in Teichm\"uller space to the central ``square" punctured torus in the moduli space of punctured tori. There are singularities in this p.d.f. arising from the topology of the moduli space.  We also consider the statistics of the distance in Teichm\"uller space to the rectangular punctured tori and the p.d.f and expected distortion of the extremal quasiconformal mappings.  \end{abstract}

\maketitle

\section{Introduction.}  In earlier work \cite{M}  we introduced a geometrically natural probability measure on the space of Riemann surfaces isometric to punctured tori.  This distribution was induced from the cross ratio distribution of the vertices (selected randomly and uniformly from the circle) of an ideal quadrilateral identified as  the fundamental domain of a rectangular punctured torus group (see \S 5 below).  This was part of a more general programme to study ``random'' discrete groups of M\"obius transformations,  see \cite{MO,MOY}.  The cross ratio distribution allowed us to calculate the basic statistics of these random punctured tori,  such as the length of the shortest geodesic and the conformal modulus.  In fact the calculation of the conformal modulus from cross ratio is a difficult classical problem initiated by Hilbert and Klein,  lies just beyond the theory of special functions,  \cite{H,KRV,K} and is intertwined with the Landau conjecture.  Thus we obtained computationally based descriptions.  However the probability distribution functions were real analytic.  In this article we consider a reverse approach.  The moduli space of punctured tori can be identified isometrically with the moduli space of lattices of $\IC$.  That space is a finite area Riemann surface - $\IS^2_{(2,3\infty)}$,  the punctured two-sphere with cone points of order two and three of hyperbolic area $\frac{\pi}{3}$.  Thus the normalised area measure will give us a uniform distribution on this moduli space.  Here we calculate statistics of the distance in Teichm\"uller space to the ``central'' square puncture torus, the punctured torus obtained by identfying the sides of an ideal hyperbolic square. It is a little surprising that there are singularities in the p.d.f. arising from the topology of the moduli space.  We then calculate the p.d.f. for distance to the rectangular punctured tori.  We obtain results of which the following is a typical example.
\begin{theorem} Let $T^{2}_{*}\in \IS^2_{(2,3\infty)}$ be a uniformly selected punctured torus.  Then there is a $K$-quasiconformal homeomorphism $f:T^{2}_{*} \to \Sigma^{2}_{*}$,  with $ \Sigma^{2}_{*}$ a rectangular punctured torus,  with $1\leq K \leq \sqrt{3}$.  Moreover the expected value $\overline{K}$  of the distortion of the extremal quasiconformal mapping is $\overline{K}=1.154\ldots$ with variance $\sigma^{2}_{K}=0.0219$.
\end{theorem}
Equation (\ref{closed}) gives closed form expressions for these numbers while the p.d.f. for the distortion of the extremal mapping is given in Theorem \ref{extremal}.  In this case the extremal quasiconformal mapping is a real analytic diffeomorphism which is Lipschitz (with constant depending on $K$) in the hyperbolic metric.  This result says that the distribution and statistics of random punctured are close to those of the (computationally approachable) rectangular punctured tori.  We refer the reader to \cite{B,G,Hb} for the theory of Teichm\"uller spaces,  though our presentation is largely self-contained and does not use any of the deep theory that has been developed.

\section{The moduli space of lattices.}
A hyperbolic punctured torus $T^{2}_{*}$  admits a covering of the form  
\begin{equation}
T^{2}_{*} \approx (\IC \setminus \Lambda(0))/\Lambda,
\end{equation} 
where  $\Lambda$ is a group of translations of $\IC$,
\begin{equation}
\Lambda =   \{z\mapsto z+ m\alpha + n \beta: m,n\in \IZ\},
\end{equation}
and $\alpha,\beta\in \IC\setminus\{0\}$,  $\alpha/\beta\not\in\IR$,  are periods.  When $\IC\setminus \Lambda(0)$ is equipped with the complete hyperbolic metric, this induces a hyperbolic metric on $T^{2}_{*} $ since $\Lambda$ acts by isometry.  Every conformal equivalence class of punctured tori can be obtained this way. 

The $SL(2,\IZ)$ action on the set of lattices allows us to normalise $\Lambda$ to assume that $\alpha=1$ and that $\beta$ lies in the fundamental domain $\Omega$ for the modular group acting on the upper half-plane, 
\begin{equation} \Omega=\Big\{z\in\IC: \Im m(z)>0, |z|\geq 1, -\frac{1}{2}< \Re e(z) \leq \frac{1}{2} \Big\}.\end{equation}
  The moduli space of lattices is the Riemann surface of area $\pi/3$, 
\[ \IH^2/PSL(2,\IZ) = \IS^{2}_{(2,3,\infty)} \]
That is the Riemann sphere with a puncture and cone points of order two and three.  When endowed with the hyperbolic metric all such $ \IS^{2}_{(2,3,\infty)} $ are isometric.  We denote the lattice with periods $1$ and $\tau\in\Omega$ by $\Lambda_\tau$.

\medskip

\scalebox{0.6}{\includegraphics[viewport=30 400 760 760]{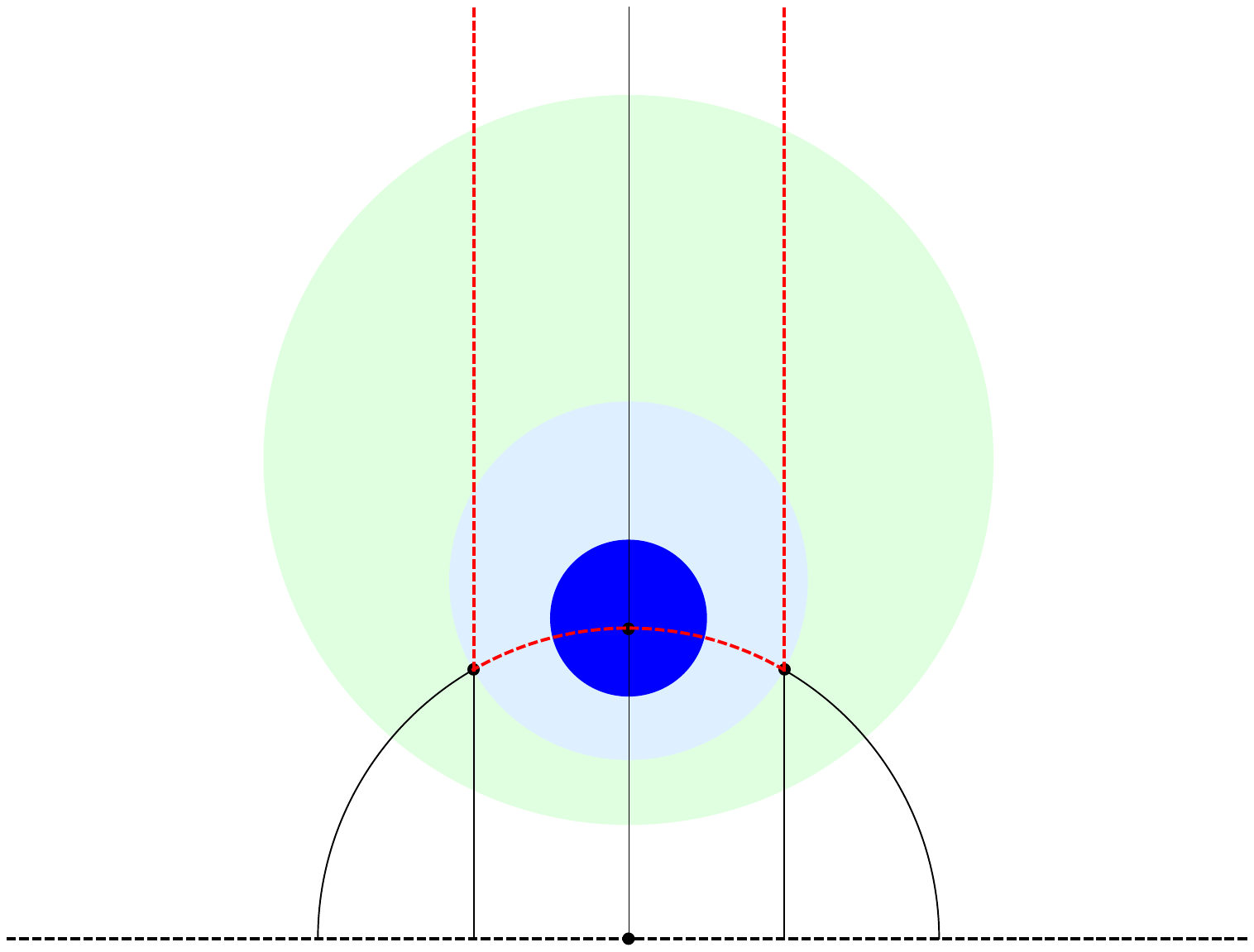}} \\

%\scalebox{0.6}{\includegraphics[viewport=30 400 760 760]{Fundamental}} \\
{\em A fundamental domain for $SL(2,\IZ)$ (hatched outline) with three disks $\ID(i,\frac{1}{4})$,  $\ID(i,\frac{1}{\sqrt{3}})$ and $\ID(i,1)$.}

\section{The Teichm\"uller metric.}
The extremal quasiconformal mapping $f:\IC\to\IC$ which commutes with the lattices $\Lambda_i$ and $\Lambda_\tau$,  that is the mapping of least distortion, is the linear mapping $(x,y)\mapsto (x+y \Re e[\tau] ,y\Im m[\tau])$.  Then $f\circ\Lambda_i=\Lambda_\tau\circ f$. In complex notation
\[ f(z) = az+b\zbar, \hskip10pt a= \big(1-i\tau\big)/2, \;\; b= \big(1+i\tau\big)/2,\]
so $f(z)=z$ if $\tau=i$.
The logarithm of the distortion of $f$ is 
\[\log K=\log \frac{|a|+|b|}{|a|-|b|}= \log \frac{1+\left| \frac{1+i\tau}{1-i\tau} \right|} {1-\left| \frac{1+i\tau}{1-i\tau} \right|} = \rho_\ID\Big(0,\frac{1+i\tau}{1-i\tau}\Big) = \rho_{\IH^2}\big(i, \tau \big)\]
where $\rho_X$ denotes the hyperbolic hyperbolic metric of curvature $-1$ on the space $X=\ID,\IH^2$. This formula relates the well know equivalence between the Teichm\"uller metric on the space of punctured tori and the hyperbolic metric of $\IS^{2}_{(2,3,\infty)}$ induced from $\IH^2$ since a lift of a quasiconformal mapping between punctured tori to $\IC\setminus \Lambda_\tau(0)$ will be quasiconformal with the same distortion and will commute with the lattices.

\section{The p.d.f. for $\rho_{\IS^{2}_{(2,3,\infty)}}\big(i^*, \tau \big)$.}
There are three different types of hyperbolic disk about the cone point of order two (call it $i^*$,  as it is the projection of $i\in\IH^2$) in $\IS^{2}_{(2,3,\infty)}$ and we will analyse each case separately.  When the radius is small,  these hyperbolic disks are isometric with a disk in the hyperbolic plane modulo a rotation of order two preserving the disk.  As soon as $\sinh(r)>\frac{1}{2}$ a hyperbolic disk about $i^*$ fails to be embedded,  and when $\sinh(r)>\frac{1}{\sqrt{3}}$ it encloses the cone point of order three as well.  The following easy lemma relates a hyperbolic disk to a Euclidean disk and is useful in making the calculations we need.
\begin{lemma}\label{disk}
\begin{equation}
\ID_{\IH^2}(i,r) = \ID(\cosh(r),\sinh(r))
\end{equation}
\end{lemma}
\noindent Since we are considering the uniform distribution on $\IS^{2}_{(2,3,\infty)}$, in each case we will need to calculate the hyperbolic area of the set $\ID_{\IS^{2}_{(2,3,\infty)}}(i^*,r)$ and this is the same as the area of $\ID_{\IH^2}(i,r)\cap \Omega$.
\subsection{$\sinh(r)\leq \frac{1}{2}$.}   In this case the unit circle includes a diameter of $\ID_{\IH^2}(i,r)$ and we easily obtain
\begin{equation}
\mbox{ h-Area}(\ID_{\IS^{2}_{(2,3,\infty)}}(i^*,r))= 2\pi\sinh^2\big(\frac{r}{2}\big)
\end{equation}
since the area of a hyperbolic disk of radius $r$ is $4\pi \sinh^2(r/2)$,  \cite{Beardon}.
\subsection{$\frac{1}{2}\leq \sinh(r)\leq \frac{1}{\sqrt{3}}$.} The hyperbolic area measure on $\IH^2$ is $|dz|^2/\Im m(z)^2$.  Using Lemma \ref{disk} we first calculate the hyperbolic area of the set $U=\ID\big(\cosh(r),\sinh(r)\big)\setminus \Omega$.  
\[ 
\mbox{\rm h-Area}(U)=2 \int_{\cosh(r)-\sqrt{\sinh^2(r)-\frac{1}{4}}}^{\cosh(r)+\sqrt{\sinh^2(r)-\frac{1}{4}}} \int_{0}^{\sqrt{\sinh^2(r)-(y-\cosh(r))^2}-\frac{1}{2}} \frac{1}{y^2} \; dx \; dy\;\;\;\;\;\; \]
\[ = 2 \cosh (r) \tan ^{-1}\Big[ \frac{4 \cosh (r) \sqrt{2 \cosh (2 r)-3}}{7-3 \cosh (2 r)}\Big]  -4 \tan ^{-1}[ \sqrt{2 \cosh (2 r)-3}]. \]
Thus $ \mbox{\rm h-Area}\big[\ID_{\IS^{2}_{(2,3,\infty)}}(i^*,r)\big] = \mbox{\rm h-Area}(\ID(\cosh(r),\sinh(r)\cap \Omega)$
\[=2\pi \sinh^2 \big[\frac{r}{2}\big] +4 \tan ^{-1}[ \sqrt{2 \cosh (2 r)-3}]   - 2 \cosh (r) \tan ^{-1}\Big[ \frac{4 \cosh (r) \sqrt{2 \cosh (2 r)-3}}{7-3 \cosh (2 r)}\Big] .
\]
\subsection{$ \sinh(r)\geq \frac{1}{\sqrt{3}}$.} The hyperbolic area of $\Omega$ is $\frac{\pi}{3}$ and so the hyperbolic area we seek now is $\frac{\pi}{3}-\mbox{h-Area}(\Omega\setminus \ID(i,r))$ and as above we calculate this to be equal to
\begin{eqnarray*}
\lefteqn{\mbox{h-Area}\big[\ID_{\IS^{2}_{(2,3,\infty)}}(i^*,r)\big]}&&\\ &=&\frac{\pi}{3}-2 \int_{0}^{1/2} \int_{\cosh(r)+\sqrt{\sinh^2(r)-x^2}}^{\infty} \;\; \frac{1}{y^2} \; dy \; dx \\
&=& \frac{\pi}{3}-\tan ^{-1}\left(\frac{1}{\sqrt{2 \cosh (2 r)-3}}\right) \\ && - \frac{1}{2} \cosh (r) \tan ^{-1}\left(\frac{4 \left(\cosh (2 r)+\cosh (r) \sqrt{2 \cosh (2 r)-3}-1\right)}{5 \cosh (2 r)-2 \cosh (4 r)-8 \sinh ^2(r) \cosh (r) \sqrt{2 \cosh (2 r)-3}+1}\right) \;\;\;\;\; \end{eqnarray*}
This last calculation provides us with the following theorem.
\begin{theorem} The cumulative distribution is  ${\bf Pr}\Big(\Big\{ d_{\IS^{2}_{(2,3,\infty)}} \, \big(\tau,i*\big) < r \Big\}\Big) =$
\begin{equation*}
=\left\{ \begin{array}{lcl}  \framebox{$0\leq r \leq \sinh^{-1}\big(\frac{1}{2}\big)$}, && \\ \\   \hskip10pt\frac{2}{3} \sinh^2\big(\frac{r}{2}\big),  \\ 
\\
 \framebox{$\sinh^{-1}\big(\frac{1}{2}\big) \leq r \leq \sinh^{-1}\big(\frac{1}{\sqrt{3}}\big)$}, && \\ \\
\hskip10pt\frac{2}{3}  \sinh ^2\left(\frac{r}{2}\right)- \frac{6}{\pi} \cosh (r) \tan ^{-1}\left[\frac{4 \cosh (r) \sqrt{2 \cosh (2 r)-3}}{7-3 \cosh (2 r)}\right]  & & \\
\;\;\;\;\hskip10pt +\frac{8}{\pi}  \tan ^{-1}\big[\sqrt{2 \cosh (2 r)-3}\big],  \\
\\
 \framebox{$r \geq \sinh^{-1}\big(\frac{1}{\sqrt{3}}\big)$}, && \\ \\
\hskip10pt1-\frac{6}{\pi}\tan ^{-1}\left[\frac{1}{\sqrt{2 \cosh (2 r)-3}}\right] && \\
\;\;\;\; \hskip10pt+\frac{1}{2} \cosh (r) \tan ^{-1}\left[\frac{4 \left(\cosh (2 r)+\cosh (r) \sqrt{2 \cosh (2 r)-3}-1\right)}{5 \cosh (2 r)-2 \cosh (4 r)-8 \sinh ^2(r) \cosh (r) \sqrt{2 \cosh (2 r)-3}+1}\right] 
\end{array} \right.
\end{equation*}
\end{theorem}

\scalebox{0.6}{\includegraphics*[viewport=50 480 760 770]{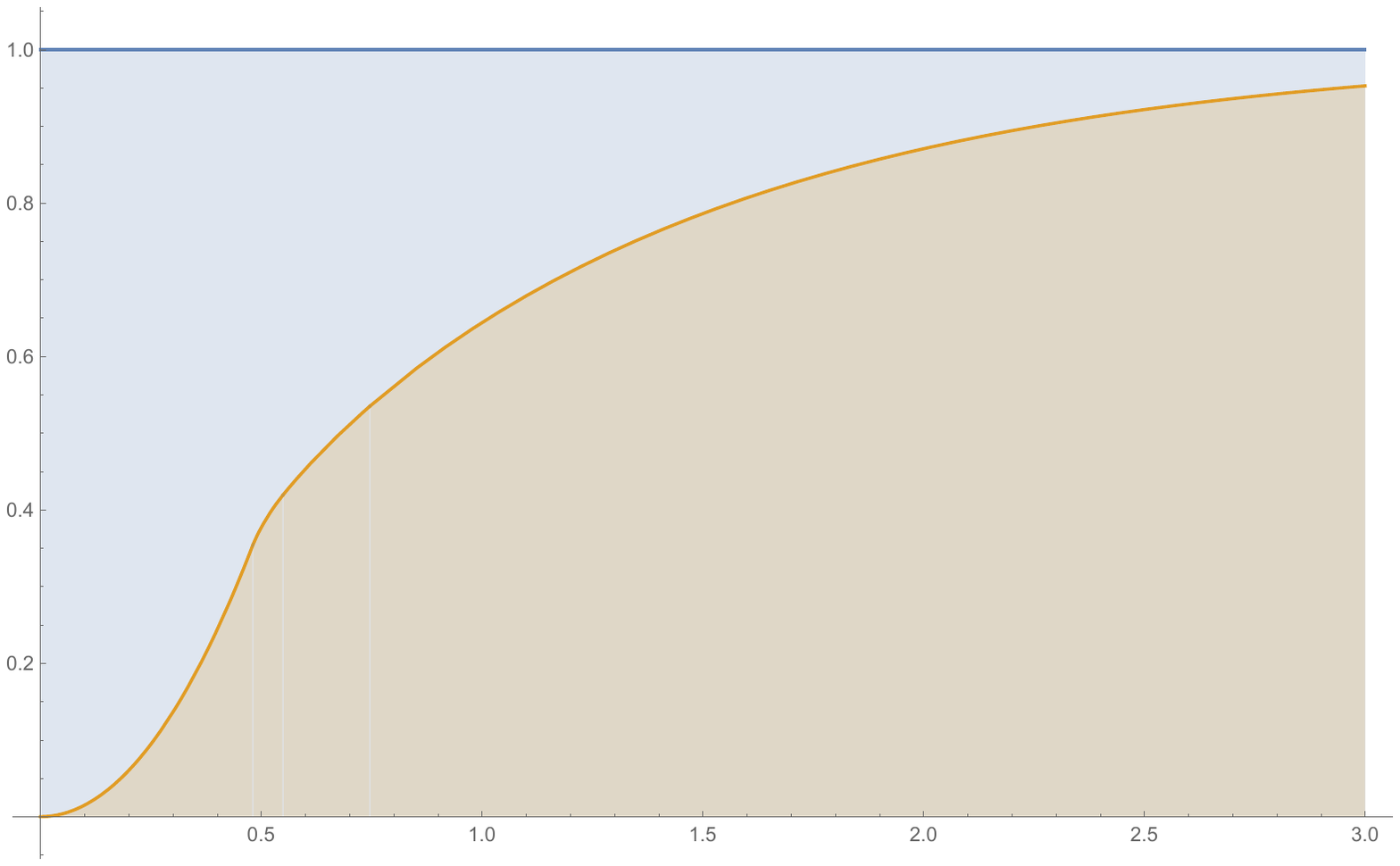}} \\
{\em The cumulative distribution of Teichm\"uller distance to the square lattice.}

\bigskip

The probability distribution function is then found from differentiating this.  Remarkably the left and right limits of the derivative at $\sinh^{-1}(1/2)$ are $3/2$ and the left and right limits of the derivative at $\sinh^{-1}(1/\sqrt{3})$ are $\frac{\sqrt{3}}{\pi} \tan ^{-1}\left(\frac{24}{7}\right)$ and so the p.d.f. is continuous.  The p.d.f. has expansion on the right of  $\sin^{-1}\big(\frac{1}{2}\big)$,
\[\frac{\pi}{2}-\frac{4 \sqrt{2}}{5^{1/4}} \sqrt{r-\sin^{-1}\big(\frac{1}{2}\big)} +O\big[r-\sin^{-1}\big(\frac{1}{2}\big)\big]. \]

\begin{theorem} The Teichm\"uller probability distribution function is the H\"older continuous function 
\begin{equation*}
=\left\{ \begin{array}{lcl}    3\sinh(r), & & \framebox{$0\leq r \leq \sinh^{-1}\big(\frac{1}{2}\big)$},   \\ 
\\
  \sinh (r) \left[3-\frac{6}{\pi} \tan ^{-1}\left[\frac{4 \cosh (r) \sqrt{2 \cosh (2 r)-3}}{7-3 \cosh (2 r)}\right]\right],   && \framebox{$\sinh^{-1}\big(\frac{1}{2}\big) \leq r \leq \sinh^{-1}\big(\frac{1}{\sqrt{3}}\big)$},   \\
\\
 \frac{3}{\pi} \sinh (r) \times &&  \framebox{$r \geq \sinh^{-1}\big(\frac{1}{\sqrt{3}}\big)$}, \\ 
\;\;\;\;  \tan ^{-1}\left[\frac{4 \left(1-\cosh (2 r)-\cosh (r) \sqrt{2 \cosh (2 r)-3}\right)}{5 \cosh (2 r)-2 \cosh (4 r)-8 \sinh ^2(r) \cosh (r) \sqrt{2 \cosh (2 r)-3}+1}\right] 
\end{array} \right.
\end{equation*}
\end{theorem}

\scalebox{0.6}{\includegraphics*[viewport=50 480 760 770]{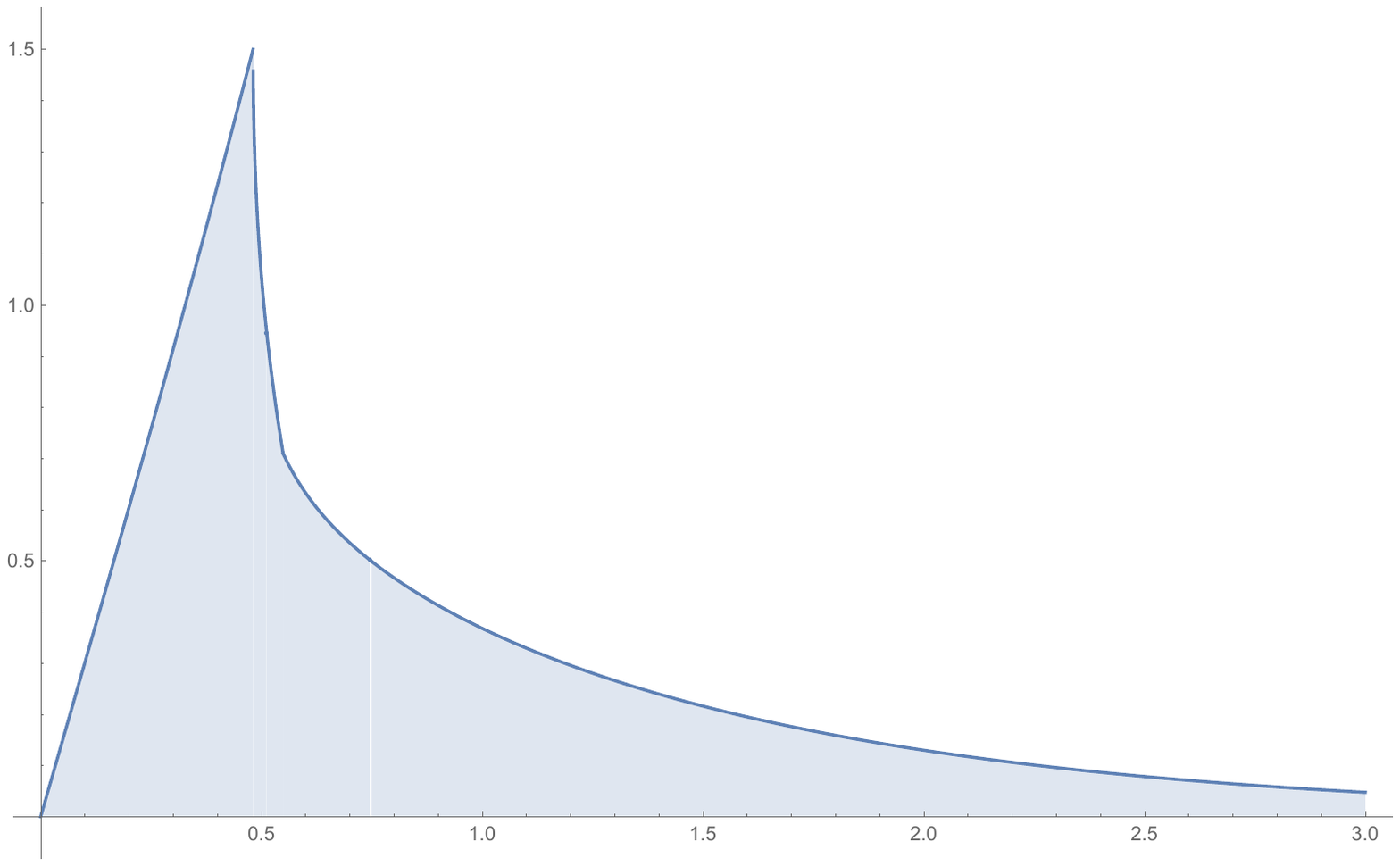}} \\
{\em The p.d.f. of Teichm\"uller distance to the square lattice.}

For large $r$ we have 
\begin{equation}\label{4.6} T(r) = \frac{3}{\pi} e^{-r} +O(e^{-5r}),  \hskip10pt r>> 1.\end{equation}
In fact this holds to 6D as soon as $r\geq 5$.  Thus all moments exist and the expected value and variance can be calculated (numerically) to be equal to
\begin{equation} E[[T]]=1.02498\ldots, \;\;\;\;  \sigma^2(T) =0.903471\ldots \end{equation}
We can interpret this result is the following way.   We may take the logarithmic transformation of this probability distribution to get the p.d.f. for $K=e^{d_{hyp}(i^*,r)}$.  Equation (\ref{4.6}) tell us that for large $K$ the p.d.f has
\[ X[K] \approx  \frac{3}{\pi K^2}+O(K^{-6}),  \]
and so has no moments.  However we could calculate the median,  but we give a few nicer examples below. 
\begin{theorem} Let $\Sigma^{2}_{*}$ denote the square punctured torus,  isometric to $(\IC \setminus \Lambda(0))/\Lambda$ with the hyperbolic metric on $ \IC \setminus \Lambda(0)$ and $\Lambda=\{z\mapsto z+m+in:m,n\in \IZ\}$.
 Let $T^{2}_{*}$ be a uniformly selected random punctured torus and let $K$ be the distortion of the extremal quasiconformal mapping $f:T^{2}_{*}\to\Sigma^{2}_{*}$.  Then
 \begin{eqnarray*}
  {\bf \rm Pr}\Big\{ K\leq \frac{1}{2}(1+\sqrt{5})\Big\} = \frac{3}{2} \left(\sqrt{5}-2\right), &&
{\bf \rm Pr}\{ K\leq 2\}   =    0.507349\ldots, 
 \end{eqnarray*}
 and  ${\bf \rm Pr}\{ K\leq 10\} =  0.904426 \ldots.$
\end{theorem}

\section{Rectangular punctured tori.} Our limited understanding of the hyperbolic metric in the complex plane with lattice points $\Lambda_\tau(0)$ removed means it is practically impossible to calculate other invariants,  for instance,  the lengths of the shortest geodesic in the punctured torus quotient space  $T^{2}_{*}(\tau)=(\IC\setminus \Lambda_\tau(0))/\Lambda_\tau$ from the modulus of a period.  However this problem is computationally feasible when $\Lambda$ is a rectangular lattice,  see \cite{Eremenko,M}.  Surprisingly general lattices are close to rectangular lattices.  We want to quantify this in a probabilistic sense.  First we give a description of the {\em rectangular} punctured tori.  

\medskip

The fundamental group of a punctured torus is free on two generators with more structure when we consider a representation as a Fuchsian group with geometric presentation $\langle a,b:[a,b]^\infty=1\rangle$,  that is the  multiplicative commutator is parabolic with $a$ and $b$ hyperbolic.  This group is represented in $PSL(2,\IC)$ by
\[ A=\pm \frac{1}{r}\, \left(
\begin{array}{cc}
 \sqrt{r^2+1} &1\\
1& \sqrt{{r^2}+1} \\
\end{array}
\right),  \hskip10pt B= \pm \frac{1}{s} \, \left(
\begin{array}{cc}
 \sqrt{{s^2}+1} & i\\
-i & \sqrt{{s^2}+1} \\
\end{array}
\right), \]
and the commutator is parabolic if and only if $rs=1$,  and then simplifies to 
\[ [A,B] = \left(
\begin{array}{cc}
 -2 i \sqrt{r^2+1} \sqrt{s^2+1}-1 & 2 i \sqrt{s^2+1}-2 \sqrt{r^2+1} \\
 -2 \sqrt{r^2+1}-2 i \sqrt{s^2+1} & 2 i \sqrt{r^2+1} \sqrt{s^2+1}-1 \\
\end{array}
\right).\]
The matrices $A$ and $B$ represent the M\"obius transformations
\begin{equation}\label{Af} f(z)=\frac{ 
 \sqrt{ {r^2}+1} \; z + 1}
 { z +\sqrt{ r^2+1}}, \hskip10pt   g(z)=\frac{  \sqrt{ {s^2}+1} z + i}{
 -i z + \sqrt{ {s^2}+1}}.
 \end{equation}
   The fixed points of $f$ are $\pm1$ and those of $g$ are $\pm i$.  Both $f$ and $g$ setwise fix the unit circle and act as isometrices of the hyperbolic disk $\ID$. The mappings $f$ and $f^{-1}$ identify the hyperbolic lines $\ell_f=\ID\cap \{z:|z + \sqrt{ {r^2}+1}|=r \}$ and $\ell_{f^{-1}} =\ID\cap \{z:|z - \sqrt{ {r^2}+1}|=r \}$, and similarly $g$ and $g^{-1}$ pair 
$\ell_g = \ID\cap\{z:|z + i\sqrt{ {s^2}+1}|=s \}$, and $\ell_{g^{-1}} = \ID\cap\{z:|z - i\sqrt{ {s^2}+1}|=s \}$.
 When $rs=1$ these four hyperbolic lines bound an ideal hyperbolic quadrilateral $Q_{r,s}$ in $\ID$ that is a fundamental polygon for the group $\langle f,g \rangle$ acting on $\ID$.  A {\em rectangular} punctured torus is any punctured torus conformally equivalent to $\ID/\langle f,g \rangle$.  These surfaces are easily identified as they come from the opposite side pairings of an ideal quadrilateral by hyperbolic elements whose axes meet at right angles.  On the surface,  there are two hyperbolic lines based at the cusp,  meeting at right angles at a finite point and which cut the surface into a quadrilateral.

 In \cite{M} we computed the probability distributions for various quantities such as the translation length of $f$ and $g$  based on the selection of a random ideal quadrilateral obtained by selecting the four vertices randomly and uniformly from $\IS$.  This distribution is not the same as that arising from the uniform distribution on $\IS^{2}_{(2,3\infty)}$,  though they do share some similarities.
 
 \begin{theorem} \label{5.7} Let the rectangular punctured torus $T^{2}_{*}$ arise from the side pairings of a random ideal quadrilateral $Q$.  Then  the shortest geodesic has p.d.f.
 \[ X_\ell = \frac{6}{\pi ^2} {\rm csch}(\ell) \Big[4 \log \cosh \frac{\ell}{2}  +2(\cosh (\ell)-1) \log \coth \frac{\ell}{2}  \Big], \hskip5pt 0< \ell \leq \log \frac{\sqrt{2}+1}{\sqrt{2}-1} , \]
with expected value $E[[\ell]]   \approx  0.984154 \ldots$.
 \end{theorem}
\begin{center}
\scalebox{0.6}{\includegraphics*[viewport=60 500 760 760]{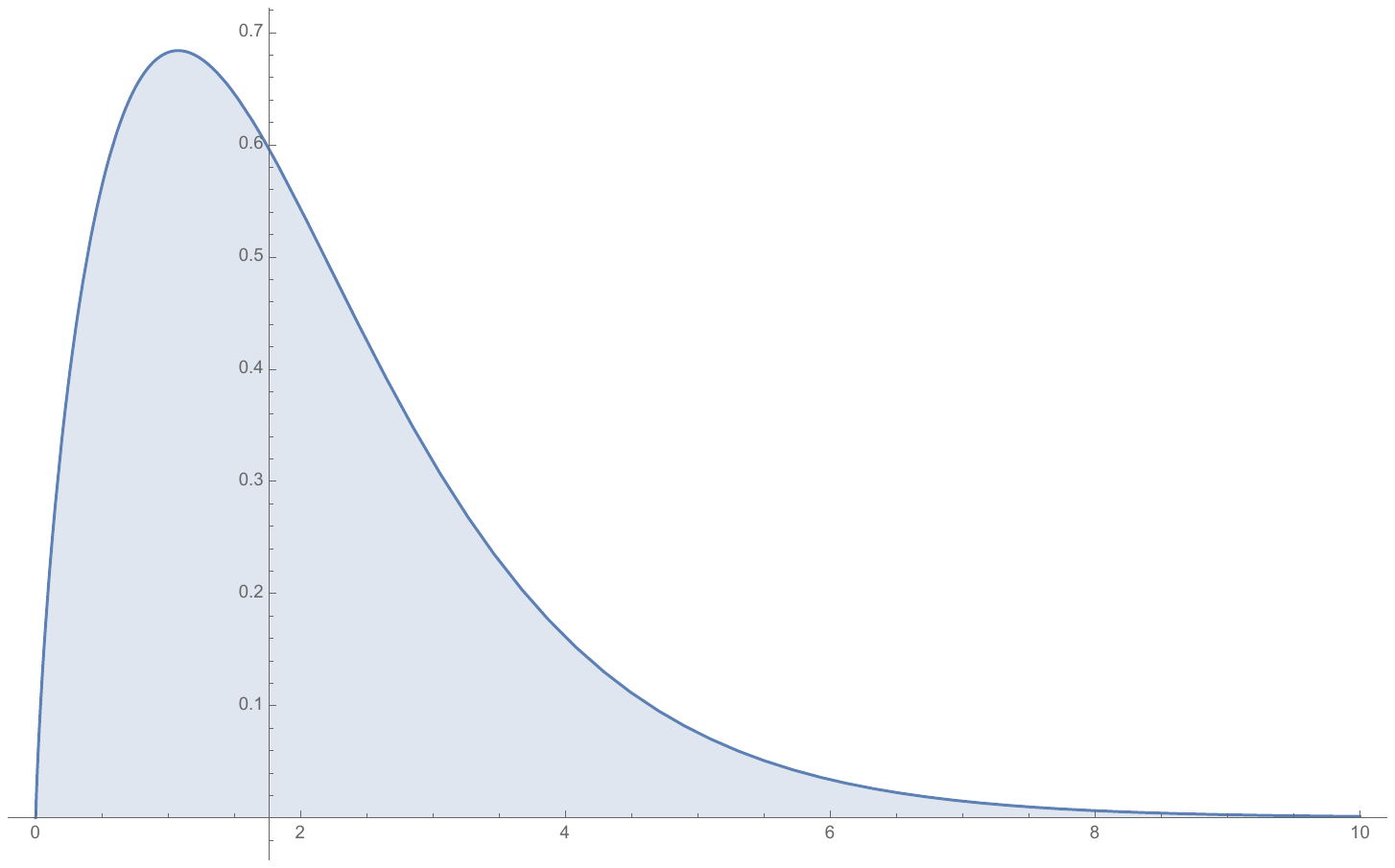}} \\
{\em Length p.d.f. : Shortest geodesic to the left of the vertical axis (at $\log \frac{\sqrt{2}+1}{\sqrt{2}-1}$),  and the length of its dual to the right.}
\end{center}
 \medskip

We therefore  now want to consider these things in regard of the uniform distribution in moduli space. We begin with some elementary hyperbolic trigonometry.

\begin{lemma}\label{lem5.1}  Let $C(r)=\{z:d_{\IH^2}(i\IR,z)\leq r\}=C_\theta$ be a cone of opening $\theta$ (the angle subtended at $0$ by $\partial C(r)$).  Then 
\begin{equation}\label{angle}
\tanh \; \frac{r}{2} = \tan \frac{\theta}{4}.
\end{equation}
\end{lemma}
\noindent{\bf Proof.} Since $\IS(0,1)\cap\IH^2$ is geodesic,  the distance $r$ is
\[ r  =  \int^{\frac{\pi}{2}}_{\frac{\pi}{2}-\frac{\theta}{2}} \; \frac{d\eta}{\sin(\eta)} =    \log \frac{1+\tan   \frac{\theta}{4}  }{1-\tan   \frac{\theta}{4}  } \] which yields (\ref{angle}) after some manipulation.  \hfill $\Box$

\medskip

The next result is immediate.

\begin{corollary}\label{cor5}  The hyperbolic distance in $\IH^2$ between $\tau$ and $|\tau|$ is
\begin{equation}\label{dist} d_{\IH^2}(\tau,|\tau|) =   2 \tanh^{-1}\Big[ \cot \Big(\frac{\arg (\tau)}{2}\Big)\Big]\end{equation}
\end{corollary}

\medskip

Next,  given $r>0$ we wish to calculate the hyperbolic area of 
\[ A(r) = \Omega \cap C(r), \] where $C(r)$ is defined in Lemma \ref{lem5.1} and $\Omega$ is the fundamental domain described earlier for the action of $SL(2,\IZ)$ on $\IH^2$.  First, given $r$ then the angle between the edge of $C(r)$ in the positive quadrant and the imaginary axis is $\eta$,  and $\tanh \frac{r}{2} = \tan \frac{\eta}{2}$ by Lemma \ref{lem5.1}.  The point of intersection of this edge with the line $\{\Re e(z)=\frac{1}{2}\}$ has modulus $a$,  where $\sin(\eta)=\frac{1}{2a}$.  Hence 
\begin{eqnarray*} \sin(\eta) &=&  \sin ( 2 \tan^{-1}(\tanh(\frac{r}{2})))  = 2 \sin (  \tan^{-1}(\tanh(\frac{r}{2})))\cos (  \tan^{-1}(\tanh(\frac{r}{2}))) \\
& = & \frac{2\tanh(\frac{r}{2})}{1+\tanh^2(\frac{r}{2})} = \tanh(r).
\end{eqnarray*}
Hence 
\[ a = \frac{1}{2\tanh(r)}.  \]
When $a\geq 1$ the hyperbolic area of the piece $\Omega\setminus \ID(0,a)$ is $A_1$ where
\begin{eqnarray*}
A_1& = & 2\int_{0}^{1/2} \int_{\sqrt{a^2-x^2}}^{\infty} \frac{dy}{y^2} \; dx = 2\int_{0}^{1/2} \frac{dx}{\sqrt{a^2-x^2}} = 2 \tan ^{-1}\left[\frac{x}{\sqrt{a^2-x^2}}\right]\Big|_{0}^{1/2}\\
& = & 2 \tan ^{-1}\left[\frac{1}{\sqrt{4a^2-1}}\right] = 2\tan^{-1}(\sinh(r)).
\end{eqnarray*}
The hyperbolic area of the piece between $\{|z|=1\}$ and $\{|z|=a\}$ is
\begin{eqnarray*}
A_2& = &  2\int^{\frac{\pi}{2}}_{\frac{\pi}{2}-\eta} \int_{1}^{a} \frac{1}{(t \sin \theta)^2} \; t \, dt\, d\theta  = (2 \log a) \int^{\frac{\pi}{2}}_{\frac{\pi}{2}-\eta} \frac{1}{\sin^2 \theta } \;  d\theta   \\
& = & 2 \log(a) \tan(\eta) = 2\sinh(r) \log\left( \frac{1}{2\tanh(r)} \right).
\end{eqnarray*}
We have now established the following theorem.
\begin{theorem}  \label{thm5} The hyperbolic area of the set $\{z\in \Omega: d_{\IH^2}(z,i\IR) \leq r \}$ is
\begin{equation}
A(r) = 2\tan^{-1}(\sinh(r)) + 2\sinh(r) \log\left( \frac{1}{2\tanh(r)} \right),
\end{equation}
if $\tanh(r)\leq \frac{1}{2}$.  Otherwise it is $\frac{\pi}{3}$.
\end{theorem}

We record the following version of Corollary \ref{cor5}.

\begin{lemma} Let $\Lambda_\tau$ be a lattice,  $\tau\in \Omega$.  Then the nearest rectangular lattice in the Teichm\"uller metric is $\Lambda_{|\tau|}$ of distance 
\begin{equation}
d_{Teich}(\Lambda_\tau,\Lambda_{|\tau|}))=d_{\IS^{2}_{(2,3,\infty)}}(\tau,|\tau|)= 2 \tanh^{-1}\Big[ \cot \Big(\frac{\arg (\tau)}{2}\Big)\Big].
\end{equation}
\end{lemma}

These calculations also give the cumulative distribution of the distance from a uniformly chosen random lattice to a rectangular lattice.  We differentiate this function (to find a remarkably simple p.d.f.) and reinterpret the result in the next theorem.

\begin{theorem}  The p.d.f. $\Sigma(r)$ for the Teichm\"uller distance of  a randomly and uniformly selected punctured torus $T^{2}_{*}$ from the moduli space,  to a rectangular punctured torus is
\begin{equation}
\Sigma(r) = \left\{ \begin{array}{ccc} \frac{6}{\pi} \; \cosh (r) \log \left(\frac{\coth (r)}{2}\right), && r\leq \tanh^{-1}(1/2), \\ 0, && r \geq \tanh^{-1}(1/2).\end{array} \right.
\end{equation}
\end{theorem}  
With a bit of work one can get the following closed form for the expected value using the Catalan number ($\approx 0.915966$).
\begin{eqnarray*} E[[\Sigma]]& = &  \frac{1}{\pi} \Big(12\; {\rm Catalan}-4 \sqrt{3} \sum _{n=0}^{\infty }\frac{(-1)^n}{3^n(2 n+1)^2}-\pi  \log (3)+12 \log \Big[\frac{\sqrt{3}+1}{2 \sqrt{2}}\Big]\Big) \\
& \approx &0.135648\ldots 
\end{eqnarray*}
and the sum is the special LerchPhi function 
\[ \Phi \left(-\frac{1}{3},2,\frac{1}{2}\right) = 4  \sum _{n=0}^{\infty }\frac{(-1)^n}{3^n(2 n+1)^2} .\]
It is also possible to get a closed form for the variance,  though it seems quite complicated.  We simply report it as 
\begin{equation}
\sigma^{2}(\Sigma) =  0.0145996\ldots
\end{equation}
Taking the Logarithmic transform of this p.d.f. gives us the p.d.f. for the distortion of the extremal quasiconformal mapping.  We record that result as follows.
\begin{theorem}\label{extremal} Let $T^{2}_{*}$ be a randomly and uniformly selected punctured torus from the moduli space $\IS^{2}_{(2,3,\infty)}$.  Then there is a $K$-quasiconformal mapping $f:T^{2}_{*}\to \Sigma^{2}_{*}$,  a rectangular punctured torus with $K\leq \sqrt{3}$.  Further the p.d.f for the distortion of the extremal quasiconformal mapping is
\[ {\bf K}= \frac{3}{\pi K^2} (K^2+1) \log \left(\frac{1}{2}\, \frac{K^2+1}{  K^2-1 }\right),  \hskip5pt  1\leq K \leq \sqrt{3}  \]
and the expected distortion is 
\begin{eqnarray}\label{closed} E[[{\bf K}]]& = & \frac{3}{\pi} \left(\frac{1}{4}  \text{Li}_2\left(\frac{1}{9}\right) -\text{Li}_2\left(\frac{1}{3}\right)+\frac{\pi ^2}{8}+\log (2)(1- \log \sqrt{3})\right) \;\;\;\; \\
& \approx & 1.15401\ldots \nonumber
\end{eqnarray}
with variance $\sigma^2({\bf K})= 0.0219564\ldots$.
\end{theorem}

\scalebox{0.7}{\includegraphics*[viewport=80 480 760 770]{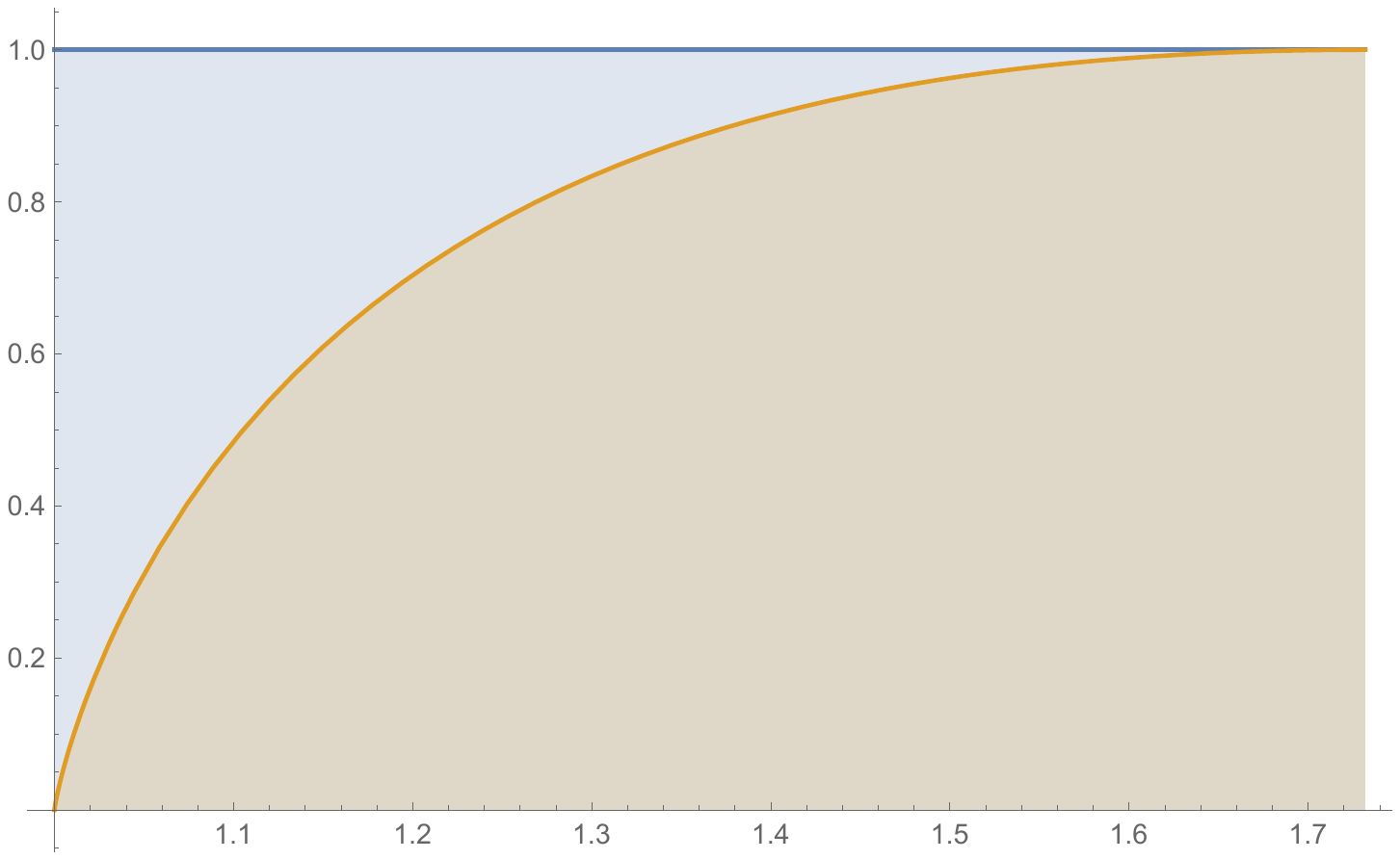}} \\
{\em The cumulative distribution of the distortion of the extremal quasiconformal mapping from a random punctured torus to a rectangular punctured torus.}

\medskip

g.j.martin@massey.ac.nz

\end{document}